\documentclass[12pt,reqno]{amsart}

\usepackage{amssymb,enumerate}

\def\thebib#1{\centerline{\normalsize \sc References}
\list
{[\arabic{enumi}]}{\settowidth\labelwidth{[#1]}\leftmargin\labelwidth
 \advance\leftmargin\labelsep
  \usecounter{enumi}}
   \def\newblock{\hskip .11em plus .33em minus .07em}
    \sloppy\clubpenalty4000\widowpenalty4000
     \sfcode`\.=1000\relax}

\DeclareMathAlphabet{\E}{U}{eus}{m}{n}     

\newcommand{\PP}{{\mathbb P}}

\newcommand{\N}{{\mathbb N}}

\newcommand{\kk}{{\Bbbk}}

\newcommand{\la}{\langle}
\newcommand{\ra}{\rangle}

\newtheorem{thm}{Theorem}[section]
\newtheorem{lemma}[thm]{Lemma}
\newtheorem{cor}[thm]{Corollary}
\newtheorem{prop}[thm]{Proposition}

\theoremstyle{definition}
\newtheorem{defn}[thm]{Definition}

\newtheorem{rmk}[thm]{Remark}

\newtheorem{example}[thm]{Example}

\newtheorem{Pf}{Proof$\!\!$}         
         
\newenvironment{pf}{\begin{Pf}}{\qed\end{Pf}}

\DeclareMathSymbol{\twoheadrightarrow}  {\mathrel}{AMSa}{"10}

\newcounter{letter}
\renewcommand{\theletter}{\rom{(}\alph{letter}\rom{)}}

\newcounter{rnum}
\renewcommand{\thernum}{\rom{(}\roman{rnum}\rom{)}}

\begin{document}


\title[Generalizing the Notion of Rank]{Generalizing the Notion of Rank to%
\\[2mm] Noncommutative Quadratic Forms}

\baselineskip15pt


\maketitle

\vspace*{0.1in}

\baselineskip15pt

\renewcommand{\thefootnote}{\fnsymbol{footnote}}
\centerline{\sc Michaela Vancliff\footnote{The first author was
supported in part by NSF grants DMS-0457022 \& DMS-0900239.\\[-3mm]}}
\centerline{Department of Mathematics, P.O.~Box 19408}
\centerline{University of Texas at Arlington,
Arlington, TX 76019-0408}
\centerline{{\sf vancliff@uta.edu}}
\centerline{{\sf www.uta.edu/math/vancliff}}

\bigskip
\centerline{and}
\bigskip

\centerline{\sc Padmini P. Veerapen\footnote{%
\begin{minipage}[t]{6.6in}%
The second author was supported in part by NSF grant DGE-0841400 as a 
Graduate Teaching Fellow in U.T.\ Arlington's GK-12 MAVS Project.%
\end{minipage}%
}}
\centerline{Department of Mathematics, P.O.~Box 19408}
\centerline{University of Texas at Arlington,
Arlington, TX 76019-0408}
\centerline{{\sf pveerapen@uta.edu}}

\setcounter{page}{0}
\thispagestyle{empty}

\bigskip
\bigskip

\begin{abstract}
\baselineskip15pt
In 2010, Cassidy and Vancliff extended the notion of a quadratic form on
$n$~generators to the noncommutative setting. In this article, we suggest a 
notion of rank for such noncommutative quadratic forms, where $n = 2$ or
$3$.  Since writing an arbitrary quadratic form as a sum of squares fails in 
this context, our methods entail rewriting an arbitrary quadratic form as 
a sum of products.  In so doing, we find analogs for $2 \times 2$ minors and 
determinant of a $3 \times 3$ matrix in this noncommutative setting. 
\end{abstract}

\baselineskip21pt


\newpage

\section*{Introduction}

Quadratic forms arise in many scientific fields, and consequently they
been studied for many decades.
Traditionally, the setting of quadratic forms has been commutative algebra
and algebraic geometry, but, in recent years, quadratic forms have played a 
role in 
noncommutative algebra via their involvement in the study of graded Clifford 
algebras (\cite{Aubry.Lemaire,L,VdB.20pts}). In this noncommutative setting, 
certain (commutative) quadratic forms associated to the graded Clifford 
algebra (GCA)
determine a quadric system $\mathfrak Q$, and the regularity of the GCA
and the degree of its generators and relations are completely determined by 
properties of 
$\mathfrak Q$. Moreover, elements of rank at most two within $\mathfrak Q$ 
determine properties of the point modules over the GCA (\cite{VVW}).

In the last few years, an algebra that is a quantized analog of a GCA was 
introduced by 
Cassidy and Vancliff in \cite{CV} and is called a graded skew Clifford
algebra (GSCA).  In this new setting, {\em noncommutative} quadratic forms 
(defined in \cite{CV}) play a role relative to GSCAs that is 
identical to that played by (commutative) quadratic forms relative to GCAs.
In particular, the regularity of the GSCA and the degree of its generators 
and relations are completely determined by properties of a certain {\em 
noncommutative} quadric system associated to the GSCA. Given $n\in\N$,
GSCAs enable the relatively-easy production of quadratic regular algebras of 
global dimension $n$.  Moreover, in \cite{CV}, many examples of GSCAs are 
given that are candidates for generic regular algebras of global 
dimension four, and, in \cite{NVZ}, it is shown that almost all quadratic 
regular algebras of global dimension three can be classified using GSCAs.

Given these recent developments, it is reasonable to attempt to extend the 
results in \cite{VVW} for GCAs to GSCAs, but, in so doing, a notion of 
{\em rank}\/ for noncommutative quadratic forms is needed.  The purpose of 
this article is to suggest such a notion of rank on the noncommutative
quadratic forms of \cite{CV}.

In Section~\ref{sec1}, we establish notation to be used throughout the 
article and outline some technical issues that motivate our approach.
Section~\ref{sec2} is devoted to the notion of rank for noncommutative
quadratic forms on two generators. Our main result of that section is
Proposition~\ref{muRk2VarsThm}, which relates the factoring of a quadratic 
form $Q$ on two generators as a perfect square to a noncommutative analog of 
the determinant of a $2 \times 2$ matrix associated to~$Q$. That result 
motivates 
our definition of rank, in Definition~\ref{muRk2VarsDefn}, of a quadratic
form on two generators.   Since our noncommutative setting depends on the
entries in a certain scalar matrix~$\mu$, our generalization of rank and
determinant are called {\em $\mu$-rank} and {\em $\mu$-determinant}, 
respectively.

The case of quadratic forms on three generators is discussed in
Section~\ref{sec3}, with our main results relating the writing of an
arbitrary quadratic form $Q$ on three generators as a sum of products to 
analogs of the $2\times 2$ minors, and determinant, of a $3 \times 3$ matrix
associated to~$Q$. In this section, our main result is 
Theorem~\ref{muRk3VarsThm},
and our definition of $\mu$-rank of a quadratic form on three generators is 
given in Definition~\ref{muRank3Vars}.

We believe it should be possible to define $\mu$-rank of a quadratic form on
$n$~generators, where $n \geq 4$, similar to our notion of $\mu$-rank in 
Definition~\ref{muRank3Vars}, where $n =3$, but we expect the methods will 
be highly computational if the $\mu$-rank is at least three. For $n$
generators, where $n \geq 4$, and $\mu$-rank at most two, the $\mu$-rank can 
be defined in terms of factoring as in Definition~\ref{rank1&2}.
Fortunately, the results in \cite{VVW} that promise 
to extend to the setting of GSCAs only entail quadratic forms of rank at 
most two; hence, the extension of those results is explored in \cite{VcV2}.

\bigskip
\bigskip


\section{Noncommutative Quadratic Forms}\label{sec1}

In this section, we set up the noncommutative setting for our quadratic
forms as defined in \cite[\S1.2]{CV}.  Our methods that are employed 
throughout the article to extend the traditional notion of rank to this 
noncommutative setting are discussed in \S\ref{issues}. 

\medskip

\subsection{Definitions}\label{defs}\hfill

Throughout the article, $\kk$~denotes an algebraically closed field such 
that char$(\kk)\neq~2$, and $M(n,\ \kk)$ denotes the vector space of 
$n \times n$ matrices with entries in $\kk$. For a graded $\kk$-algebra~$B$, 
the span of the homogeneous elements in $B$ of degree $i$ will be denoted 
$B_i$, and the notation $T(V)$ will denote the tensor algebra on the 
vector space~$V$.  If $C$ is any ring or vector space, then $C^\times$ will 
denote the nonzero elements in $C$. We use $R$ to denote the polynomial
ring on degree-one generators $x_1, \ldots , x_n$.

\medskip

For $\{i,\ j\} \subset \{1, \ldots , n\}$, let $\mu_{ij} \in \kk^{\times}$ 
satisfy the property that $\mu_{ij}\mu_{ji} = 1$ for all $i \ne j$. We write 
$\mu = (\mu_{ij}) \in M(n,\ \kk)$.  As in \cite{CV}, we write $S$ for the 
quadratic $\kk$-algebra on generators $z_1, \ldots, z_n$ with defining 
relations $z_j z_i = \mu_{ij} z_i z_j$ for all $i$, $j = 1, \ldots, n$, where 
$\mu_{ii} = 1$ for all $i$. 

\begin{defn}\cite[\S1.2]{CV}
\begin{enumerate}
\item[(a)]
With $\mu$ and $S$ as above, a quadratic form $Q$ is any element of $S_2$.
\item[(b)]
A matrix $M \in M(n,\ \kk)$ is called $\mu$-symmetric if $M_{ij} = 
\mu_{ij}M_{ji}$ for all $i$, $j = 1, \ldots , n$. 
\end{enumerate}
\end{defn}
\noindent 
We write $M^{\mu}(n,\ \kk)$ for the set of $\mu$-symmetric matrices in 
$M(n ,\ \kk)$.  Clearly, if $\mu_{ij} = 1$ for all $i,\ j$, then 
$M^\mu (n ,\ \kk)$ consists of all symmetric matrices.  

Henceforth, we assume that $\mu_{ii} = 1$ for all $i$.

As was shown in \cite[\S1.2]{CV}, the one-to-one correspondence between 
commutative quadratic forms and symmetric matrices has a counterpart in
our setting, with the one-to-one correspondence being between
noncommutative quadratic forms and $\mu$-symmetric matrices.  This
correspondence is given as follows. For $M \in M^{\mu}(n,\ \kk)$, let 
$\tilde{Q} = z^T M z \in T(S_1)_2$, where $z = (z_1, \ldots, z_n)^T$, and  
write $Q$ for the image in $S$ of $\tilde{Q}$; the element $Q$ is the
quadratic form corresponding to $M$. Conversely, if 
$Q = \sum_{i \leq j} \alpha_{ij} z_i z_j \in S$, where $\alpha_{ij} \in
\kk$ for all $i,\ j$, is a quadratic form, then the matrix $(M_{ij})$,
where $M_{kk} = \alpha_{kk}$, $M_{ij} = 2^{-1}\alpha_{ij}$ and 
$M_{ji} = 2^{-1}\mu_{ji} \alpha_{ij}$ for all $i,\ j,\ k$ where $i < j$, 
is the $\mu$-symmetric matrix corresponding to $Q$.

\begin{rmk}\label{twistrmk}
As was shown in \cite{ATV1}, if the point modules of $S$ are parametrized by
$\PP^{n-1}$, then $S$ is a twist (in the sense of \cite[\S8]{ATV2}) of the 
polynomial ring by a graded degree-zero automorphism $\tau \in$ Aut$(R)$
(see Definition~\ref{twistdef} below).  This case occurs if and only if 
$\mu_{ik} = \mu_{ij}\mu_{jk}$ for all $i,\ j,\ k$. This is the situation
throughout Section~\ref{sec2}, since there the assumption that $n=2$ causes 
the point modules of $S$ to be parametrized by $\PP^1$.
\end{rmk}

\medskip

\subsection{Technical Issues}\label{issues}\hfill

Recall that, since $\kk$ is algebraically closed, every (commutative)
quadratic form $q \in R^\times$ can be written as $\sum_{i=1}^m y_i^2$, where 
$y_1, \ldots , y_m \in R_1$ are linearly independent and $m \in \N$ is
unique; in this setting, the rank of $q$ is defined to be $m$.
However, if $Q \in S_2$ is a noncommutative quadratic form, then a direct 
generalization using a sum of squares leads to problems that are
demonstrated by the following example.
\begin{example}
Suppose $n = 2$, $\mu_{12} = -1$ and $Q = z_1^2 + 2bz_1z_2 + cz_2^2$,
where $b$, $c \in \kk$. If $b \ne 0$, then $Q \ne \sum_{i=1}^m X_i^2$ for any 
$m \in \N$, where $X_i \in S_1$ for all $i$.
Moreover, if $b = 0$, then 
$Q = z_1^2 + c z_2^2 = (z_1 + \alpha z_2)^2$, where $\alpha \in \kk$,
$\alpha^2 = c$. Hence, if $b \neq 0$, then a sum of squares is not possible; 
whereas if $b=0$, then a sum of square terms is possible but the number of 
such terms is not unique.
\end{example}

Instead, for $n \leq 3$, we will model our notion of rank on the 
following facts concerning the rank of a (commutative) quadratic form 
$q \in R_2$: 
\begin{enumerate}
\item[(a)] rank$(q) = 0$ if and only if $q = 0$;
\item[(b)] rank$(q) = 1$ if and only if $q =  X^2$ for some 
            $X \in R_1^\times$;
\item[(c)] rank$(q) = 2$ if and only if $q =  XY$ for some linearly
                 independent  $ X,\ Y \in R_1$;
\item[(d)] rank$(q) = 3$ if and only if $q =  XY + Z^2$ for some linearly
                 independent $ X,\ Y,\ Z \in R_1$.
\end{enumerate}
However, in our noncommutative setting, it is possible that $Q \in S_2$ might
factor as both a perfect square and also as a product of linearly independent 
elements. This issue is highlighted in the next example.

\begin{example}\label{NUF2Vars}
If $Q = z_1^2 + 6z_1z_2 + 4z_2^2 \in S_2$, with $\mu_{12} = 2$, then
$Q = (z_1 + 2 z_2)^2 = (z_1 + z_2)(z_1 + 4z_2)$.
\end{example}


\bigskip
\bigskip

\section{Rank of Quadratic Forms on Two Generators}\label{sec2}

In this section, we consider noncommutative quadratic forms on two generators
as defined in Section~\ref{defs}, and introduce a notion of rank, called
$\mu$-rank, in Definition~\ref{muRk2VarsDefn}, on such quadratic forms that 
extends the notion of rank of a commutative quadratic form. We also
introduce in Definition~\ref{D(M)} an analog of the determinant of a 
$2 \times 2$ matrix.

Throughout this section we suppose $n = 2$.
As mentioned in Remark~\ref{twistrmk}, we use the notion of {\em twist} in 
this section, which is defined as follows.

\begin{defn}\label{twistdef}\cite[\S8]{ATV2} \
Let $B = \bigoplus_{m\geq 0} B_m$ be a quadratic algebra and let $\phi$
be a graded degree-zero automorphism of $B$.  The twist $B^\phi$ of $B$ by
$\phi$ is the vector space $\bigoplus_{m \geq 0} B_m$ with a new
multiplication $*$ defined as follows: if $x, y \in B_1$, then
$x * y = x \phi(y)$, where the right-hand side is computed using the
original multiplication in $B$.
\end{defn}

By Remark~\ref{twistrmk}, in this section, the algebra $S$ is a twist of the 
polynomial ring $R$ by a graded, degree-zero automorphism $\tau \in$ 
Aut$(R)$.  In this section, we denote multiplication in~$S$ by $*$ and the 
action of $\tau$ by $r^\tau = \tau(r)$ for all $r\in R$.
By \cite[Lemma 5.6]{N}, we may choose $\tau$ to be given by  
\[
\tau(z_1) = \mu_{12} z_1 \quad \text{and} \quad \tau(z_2) = z_2. \tag{$*$}
\]

\begin{lemma}
Suppose $n = 2$. If $Q \in S_2$ is a quadratic form, then $Q$ factors in 
at most two distinct ways.
\end{lemma}
\begin{pf}
Suppose $Q = r_1 *r_2 = r_3 * r_4 = r_5 * r_6$ in $S$, where $r_i \in S_1$ 
for all $i$. Using $\tau$ given above in \thetag{$*$}, it follows that 
$Q = r_1r_2^{\tau} = r_3r_4^{\tau} = r_5r_6^{\tau}$ in $R$.
However, in $R$, the element $Q$ factors in at most two distinct ways, so, 
without loss
of generality, we may assume $r_5 \in \kk^{\times} r_3$ and $r_6 \in 
\kk^{\times} r_4$. Hence, in $S$, $Q$ factors in at most two ways. 
\end{pf}

For the rest of this section, we will be concerned with a quadratic form 
$a z_1 *z_1 + 2b z_1 * z_2 + c z_2*z_2 \in S_2$, where $a,\ b,\ c \in \kk$. 
As explained in \S\ref{defs}, to such a quadratic form is associated a
$\mu$-symmetric matrix $M = 
\begin{bmatrix} a & b \\ \mu_{21} b & c \end{bmatrix}$. It will be
useful to use an analog of the determinant function on $M$ in the next result.

\begin{defn}\label{D(M)}
Let $D: M^{\mu}(2,\ \kk) \to \kk$ be given by 
\[
D(M) = 4b^2 - (1+\mu_{12})^2 ac, \quad \text{where}\quad
M = \begin{bmatrix} a & b \\ \mu_{21} b & c \end{bmatrix};
\] 
we call $D(M)$ the $\mu$-determinant of $M$.
\end{defn}
\noindent
We remark that if $S = R$, that is, if $\mu_{12} = 1$, then 
$D(M) = -4 \det (M)$.

\begin{prop}\label{muRk2VarsThm}
Let $Q = a z_1 *z_1 + 2b z_1 * z_2 + c z_2*z_2 \in S_2^\times$,
where $a,\ b,\ c \in \kk$, be a quadratic form with associated
$\mu$-symmetric matrix $M \in M^{\mu}(2,\ \kk)$.
\begin{enumerate}
\item[{\rm (a)}] There exists $L_1, L_2 \in S_1$ such that $Q = L_1 *L_2$ 
                in $S$. 
\item[{\rm (b)}] There exists $L \in S_1$ such that $Q = L*L$ in $S$
      if and only if $D(M) = 0$.
\item[{\rm (c)}] The element $Q$ factors uniquely, up to a nonzero
                 scalar multiple, in $S$ if and only if $b^2 = \mu_{12} ac$.
\end{enumerate}
\end{prop}
\begin{pf}
Viewing $Q \in R$, we have $Q = a \mu_{12} z_1^2 + 2b z_1 z_2 + c z_2^2$.

(a) Since $Q$ factors in $R$, we have $Q = r_1 r_2$, 
where $r_i\in R_1 = S_1$ for all $i$.
Thus, in $S$, $Q = r_1 * \tau^{-1}(r_2)$, which proves (a). 

(b) If $Q = r * r$ in $S$, for some $r \in S_1$, then 
\[
Q = r r^\tau = \mu_{12} \alpha_1^2 z_1^2 + (1 + \mu_{12}) \alpha_1
\alpha_2 z_1 z_2 + \alpha_2^2 z_2^2
\]
in $R$, where $r = \alpha_1 z_1 + \alpha_2 z_2$ for some $\alpha_1,\ 
\alpha_2 \in \kk$. Comparing coefficients, it follows that this situation
occurs if and only if $2b = (1 + \mu_{12}) \alpha_1 \alpha_2$,
where $\alpha_1^2 = a$ and $\alpha_2^2 = c$. Hence,  $Q = r*r$ for some
$r \in S_1$  implies that $D(M) = 0$. Conversely, if $D(M) = 0$, then 
$2 b = (1 + \mu_{12}) \beta$, where $\beta \in \kk$ and $\beta^2 = ac$.
If also $ac = 0$, then (b) follows; whereas if $ac \neq 0$, then we may
choose $\alpha_1,\ \alpha_2 \in \kk$ such that $\alpha_1^2 = a$ and
$\alpha_2 = \beta/ \alpha_1$, which implies that $Q = r*r$ in $S$, where 
$r = \alpha_1 z_1 + \alpha_2 z_2$.

(c) A quadratic form factors uniquely in $S$ if and only if it factors 
uniquely in $R$, and the latter occurs if and only if the discriminant
is zero. Since the discriminant of $a \mu_{12} z_1^2 + 2b z_1 z_2 + c z_2^2
\in R_2$ belongs to $\kk^\times(b^2 - \mu_{12}ac)$, the result follows.
\end{pf}

\begin{cor}
Let $Q$ be as in Proposition~\ref{muRk2VarsThm}.  
\begin{enumerate}
\item[{\rm (a)}]
Suppose $Q$ does not factor uniquely. If $ac = 0$, then $Q \in \la z_i \ra$ 
for some $i \in \{1,\ 2\}$%
$;$ whereas if $ac \neq 0$, then 
\[
Q = \left(z_1 + \frac{c z_2}{b + H}\right)*
     \left(a\, z_1 + \left[b + H \right]z_2\, \right),
\]
where $H^2 = b^2 - \mu_{12}ac$.
\item[{\rm (b)}]
Suppose $Q$ factors uniquely, up to a nonzero scalar multiple, in $S$.
If $b = 0$, then $Q \in \kk^\times z_i^2$ for some $i \in \{1,\ 2\}$%
$;$ whereas if $b \neq 0$, then 
\[
Q = b^{-1} (b z_1 + c z_2)*(a z_1 + b z_2).
\]
\end{enumerate}
\end{cor}
\begin{pf}\hfill

(a) If $ac = 0$, the result in (a) clearly holds. 
If $ac \neq 0$, we may write 
$Q = a^{-1}(a z_1 + \alpha z_2)*(a z_1 + \beta z_2)$, 
where $\alpha$, $\beta \in \kk^\times$. Comparing coefficients, we find
$ac = \alpha \beta$ and $2 b = \beta + \mu_{12} \alpha$.
Solving for $\beta$ yields $\beta = b + H$, where $H^2 = b^2 - \mu_{12}
ac$. Since $\alpha = ac/(b+H)$, part (a) follows.

(b) By Proposition~\ref{muRk2VarsThm}(c), $b^2 = \mu_{12}ac$. 
Thus, if $b = 0$, the result in (b) clearly holds. 
If $b \neq 0$, then $ac \neq 0$, so part (a) applies with $H = 0$.
\end{pf}

Proposition~\ref{muRk2VarsThm} suggests the following generalization of 
the rank of a quadratic form on two generators.

\begin{defn}\label{muRk2VarsDefn}
Let $Q = a z_1 *z_1 + 2b z_1 * z_2 + c z_2*z_2 \in S_2$, where 
$a$, $b$, $c\in \kk$, let $M\in M^{\mu}(2,\ \kk)$ be the $\mu$-symmetric 
matrix associated to $Q$ and let 
$D: M^{\mu}(2,\ \kk) \to \kk$ be defined as in Definition~\ref{D(M)}. 
If $n=2$, we define $\mu$-rank $: S_2 \to \N$ as follows:
\begin{enumerate}[(a)]
  \item if $Q = 0$, we define $\mu$-rank$(Q) = 0$;
  \item if $Q \neq 0$ and $D(M) = 0$, we define $\mu$-rank$(Q) = 1$;
  \item if $D(M) \ne 0$, we define $\mu$-rank$(Q) = 2$.
\end{enumerate}
\end{defn}
\begin{example}
If $Q$ is the quadratic form in Example~\ref{NUF2Vars}, then 
$\mu$-rank$(Q) = 1$.
\end{example}
\begin{cor}
Let $n = 2$.  If $Q \in S_2^\times$, then $\mu$-rank$(Q)=1$ if and only if 
$Q = L*L$ for some $L \in S_1^\times$. 
\end{cor}
\begin{pf}
Combine Definition~\ref{muRk2VarsDefn} and Proposition~\ref{muRk2VarsThm}(b).
\end{pf}


\bigskip
\bigskip

\section{Rank of Quadratic Forms on Three Generators}\label{sec3}

In this section, we explore further the notion of rank on noncommutative 
quadratic forms, and extend the results of the previous section concerning 
$\mu$-rank of quadratic forms on two generators to quadratic forms on three 
generators. Our main result of this section is Theorem~\ref{muRk3VarsThm},
which uses analogs of the determinant and minors of a $3 \times 3$ matrix
to describe factoring properties of a quadratic form. Our definition of 
$\mu$-rank of a noncommutative quadratic form on three generators is given in
Definition~\ref{muRank3Vars}.

Since $n=3$ throughout this section, the methods of Section~\ref{sec2} cannot
be employed directly since the algebra $S$, where $n \geq 3$, need not be a 
twist of a polynomial ring. In particular, we henceforth use juxtaposition
to denote the multiplication in $S$.

\begin{prop}\label{L1L2L3^2}
If $Q = az_1^2 + bz_2^2 + cz_3^2 + 2dz_1z_2 + 2ez_1z_3 + 2fz_2z_3 \in S_2$, 
where $a, \ldots , f \in \kk$, is a quadratic form, then 
$Q = L_1 L_2 + L_3^2$ for some $L_1$, $L_2$, $L_3 \in S_1$.
\end{prop}
\begin{pf}
If $a = b = c = e = 0$, then the result clearly holds. Moreover, if 
$a = b = c = 0 \neq e$, then 
\[ Q = 2(z_1 + \alpha z_2)(d z_2 + e z_3) - 2\alpha d z_2^2, \]
where $\alpha \in \kk$ and $\alpha e = f$. Hence, by symmetry, it suffices 
to prove the result in the case $a \neq 0$. Thus, for simplicity, we 
henceforth assume that $a = 1$.

If $\mu_{12} \ne -1 \ne \mu_{13}$, then 
\[ Q = Q' + \left(z_1 + \frac{2d}{1+\mu_{12}} z_2 + 
\frac{2e}{1+\mu_{13}} z_3\right)^2
,\]
where $Q' \in \kk z_2^2 + \kk z_3^2 + \kk z_2z_3$.  Applying
Proposition~\ref{muRk2VarsThm}(a) to $Q'$ implies the result in this case.

Suppose $\mu_{12} = -1 \ne \mu_{13}$. If $c \neq 0$ or $e \neq 0$, then
there exists $\delta \in \kk$ such that $\delta^2 = c$ and 
$2e \neq (1 + \mu_{13})\delta$. In this case, 
\[
Q = (z_1 + \gamma z_2 + \delta z_3)^2 + (z_1 + \alpha z_2)(2 d z_2 +
\beta z_3),
\]
where  $\alpha, \ldots , \delta \in \kk$ satisfy 
\[
\delta^2 = c,\quad 
\beta = 2 e - (1 + \mu_{13})\delta \neq 0,\quad
\gamma^2 = b - 2 d \alpha \quad\text{and}\quad
(1 + \mu_{23}) \gamma \delta + \alpha \beta = 2f.\] 
However, if $c = 0 = e$, then 
$
Q = (z_1 + \epsilon z_2)^2 - 2 z_2 (d z_1 - f z_3)$,
where $\epsilon \in \kk$, $\epsilon^2 = b$. Similarly, if $\mu_{12} 
\neq -1 = \mu_{13}$.

It remains to consider $\mu_{12} = -1 = \mu_{13}$. If $e \neq 0$, then
there exist solutions $\alpha$, $\beta$, $\gamma \in \kk$ to the equations 
\[ 
\alpha^2 + 2 d \gamma = b, \quad 
\beta^2 = c \quad \text{and} \quad
(1 + \mu_{23}) \alpha \beta + 2 e \gamma = 2f,
\]
so that 
\[
Q = (z_1 + \alpha z_2 + \beta z_3)^2 + 2 (z_1 + \gamma z_2)(d z_2 + e z_3).
\]
On the other hand, if $e = 0$, then 
$
Q = (z_1 + \delta z_3)^2 + (2 d z_1 + b z_2 + 2 \mu_{32} f z_3) z_2,
$
where $\delta \in \kk$, $\delta^2 = c$.
\end{pf}

In order to generalize Proposition~\ref{muRk2VarsThm} and 
Definition~\ref{muRk2VarsDefn} to the three-generator case, we introduce 
analogs of the determinant and $2 \times 2$ minors of a $3 \times 3$ matrix.

\begin{defn}\label{defnDi}
Let $M =
\left[\begin{smallmatrix}
    a & d & e \\[1mm]
    \mu_{21}d & b & f \\[1mm]
    \mu_{31}e & \mu_{32}f & c \end{smallmatrix}
\right] \in M^{\mu}(3, \kk)$ and, for $1 \le i \le 8$,
define the functions $D_i: M^{\mu}(3, \kk) \to \kk$ by\\[-3mm]
\begin{gather*}
\begin{array}{ll}
D_1(M) = 4d^2 - (1 + \mu_{12})^2 ab, \qquad &
D_4(M) = 2(1 + \mu_{23})de - (1 + \mu_{12})(1 + \mu_{13})af,\\[2mm]
D_2(M) = 4e^2 - (1 + \mu_{13})^2 ac, &
D_5(M) = 2(1 + \mu_{12})ef - (1 + \mu_{13})(1 + \mu_{23})cd,\\[2mm]
D_3(M) = 4f^2 - (1 + \mu_{23})^2 bc, &
D_6(M) = 2(1 + \mu_{13})df - (1 + \mu_{12})(1 + \mu_{23})be,
\end{array}\\[2mm]
\begin{array}{l}
D_7(M) = (\mu_{23}cd^2 - 2def + be^2)
            (\mu_{13}\mu_{21}cd^2 - 2def + \mu_{12}\mu_{23}\mu_{31}be^2),
	         \\[3mm]
D_8(M) = \mu_{21}(d + X)(e - Y) + \mu_{23}\mu_{31} (d-X)(e+Y)-2af,
\end{array}
\end{gather*}
\quad\\[-3mm]
where $X^2 = d^2 - \mu_{12} ab$ \ and \ $Y^2 = e^2 - \mu_{13} ac$.
We call $D_1, \ldots , D_6$ the $2\times 2$ $\mu$-minors of $M$. The
functions $D_7$ and $D_8$ will play a role analogous to that of the
determinant of $M$ and so could be called the $\mu$-determinants of $M$,
even though $D_8$ is not a polynomial in the entries of $M$. (Attempting
to convert $D_8$ to a polynomial leads to unwieldy polynomials such as
the one given after Theorem~\ref{muRk3VarsThm}.)
\end{defn}

\begin{thm}\label{muRk3VarsThm}
Let $Q = az_1^2 + bz_2^2 + cz_3^2 + 2dz_1z_2 + 2ez_1z_3 + 2fz_2z_3 \in
S_2$, where $a, \ldots, f \in \kk$, and let $M \in M^{\mu}(3,\ \kk)$ be the 
$\mu$-symmetric matrix associated to $Q$. 
\begin{enumerate}[(a)]
\item[{\rm (a)}] 
There exists $L \in S_1$ such that $Q = L^2$
if and only if $D_i(M) = 0$ for all $i = 1, \ldots , 6$.
\item[{\rm (b)}] 
\begin{enumerate}
\item[{\rm (i)}] 
If $a = 0$, then there exists $L_1$, $L_2 \in S_1$ such that
$Q = L_1 L_2$ if and only if $D_7(M) = 0$;
\item[{\rm (ii)}] 
if $a \neq 0$, then there exists $L_1$, $L_2 \in S_1$ such that
$Q = L_1 L_2$ if and only if $D_8(M) = 0$ for some $X$ and $Y$ satisfying
$X^2 = d^2 - \mu_{12} ab$ \ and \ $Y^2 = e^2 - \mu_{13} ac$.
\end{enumerate}
\end{enumerate}
\end{thm}
\begin{pf}
By Proposition~\ref{L1L2L3^2},  $Q = L_1 L_2 + L_3^2$ for some 
$L_1$, $L_2$, $L_3 \in S_1$. 

(a) \  
Suppose there exist $\alpha_1$, $\alpha_2$, $\alpha_3 \in \kk$ such that 
$Q = (\alpha_1 z_1 + \alpha_2 z_2  + \alpha_3 z_3)^2$. Comparing
coefficients, it follows that\\[-4mm]
\begin{center}
\begin{tabular}{rlrl}
(i) & $2d = (1 + \mu_{12})\alpha_1\alpha_2$, \qquad \qquad \qquad & 
              (iv)& $a = \alpha_1^2$,\\[3mm]
(ii) & $2e = (1 + \mu_{13})\alpha_1\alpha_3$, &(v)& $b = \alpha_2^2$,\\[3mm]
(iii) & $2f = (1 + \mu_{23})\alpha_2\alpha_3$, &(vi)& $c = \alpha_3^2$,
\end{tabular}
\end{center}
\quad\\[-2mm]
so $D_i(M) = 0$ for $i = 1$, $2$, $3$. Moreover, from equations (i)-(iv),
we have
\[
\begin{array}{rl}
4de (1 + \mu_{23})& = (2d)(2e)(1 + \mu_{23})\\[3mm] 
                  & = (1 + \mu_{12})(1 + \mu_{13}) (1 + \mu_{23}) 
                        \alpha_1^2\alpha_2\alpha_3\\[3mm]
                  & = (1 + \mu_{12})(1 + \mu_{13}) 2 a f,
\end{array}
\]
so $D_4(M) = 0$. By symmetry, $D_i(M) = 0$ for $i = 5$, $6$.

Conversely, suppose that $D_i(M) = 0$ for all $i = 1, \ldots , 6$. If 
$a = 0$, then $d = 0 = e$, since $D_1(M) = 0 = D_2(M)$. In this case, 
$Q \in \kk z_2^2 + \kk z_3^2 + \kk z_2 z_3$, so 
Proposition~\ref{muRk2VarsThm}(b) applies to~$Q$ (since $D_3(M) = 0$), and
so $Q = L^2$, where $L \in S_1$. Thus, to complete the proof of (a), 
we may assume $a \neq 0$.

Since $D_i(M) = 0$ for $i = 1$, $2$, $3$, there exist $w_1$, $w_2$, $w_3 
\in \kk$ such that  
\[
2d = (1 + \mu_{12}) w_1,\qquad 
2e= (1 + \mu_{13}) w_2,\qquad 
2f= (1 + \mu_{23}) w_3, \tag{vii}
\]
where 
$w_1^2 = ab$, $w_2^2 = ac$, $w_3^2 = bc$.
Since $a \neq 0$, let $Q' = a^{-1}(a z_1 + w_1 z_2 + w_2 z_3)^2
\in S_2$. By (vii), it follows that
\[
Q' = a z_1 ^2 + b z_2^2 + c z_3^2 + 2d z_1 z_2 + 2 e z_1 z_3 + 
            a^{-1} (1 + \mu_{23}) w_1 w_2 z_2 z_3. 
\]
If $ (1 + \mu_{23}) bc = 0$, then $Q' = Q$
and (a) follows. If $\mu_{12} = -1$, then $w_1$ may be chosen so that 
$Q' = Q$; similarly for $w_2$ if $\mu_{13} = -1$. Hence, we may assume 
\[
(1 + \mu_{12}) (1 + \mu_{13}) (1 + \mu_{23}) bc \neq 0. \tag{viii}
\]
Moreover,\\[-3mm] 
\[
\begin{array}{rcll}
(1 + \mu_{12}) (1 + \mu_{13}) (1 + \mu_{23}) w_1 w_2 & =&4 de (1 + \mu_{23}),
& \text{using \thetag{vii}}\\[2mm]
   &=&2 (1 + \mu_{12}) (1 + \mu_{13}) af, & \text{as } D_4(M) = 0
                 \\[2mm]
   &=&(1 + \mu_{12}) (1 + \mu_{13}) (1 + \mu_{23}) a w_3, \quad & \text{using 
               \thetag{vii}}.
\end{array}
\]
\quad\\[-2mm]
Thus, since \thetag{viii} holds, $w_1 w_2 = a w_3$, from which it follows 
that $Q' = Q$, which completes the proof of (a).

%
(b)(i) \  
Suppose $a = 0$. If also $d = 0$, then, by 
Proposition~\ref{muRk2VarsThm}(a),  $Q$ factors if and only if $be = 0$, and
the latter holds if and only if $D_7(M) = 0$. Since a similar argument
applies if instead $a = 0 = e$, we may assume $de \neq 0$. Let 
$Q_1,\ Q_2 \in S_2$ be given by\\[-2mm] 
\[
\begin{array}{rcl}
Q_1 &=& 2 [z_1 + (2d)^{-1}b z_2 + (2e)^{-1}c z_3][d z_2 + e z_3]\\[3mm]
    &=& b z_2^2 + c z_3^2 + 2 d z_1 z_2 + 2 e z_1 z_3 + 
         ( be d^{-1} + cd \mu_{23} e^{-1} ) z_2 z_3,\\[5mm]
Q_2 &=& 2 [d \mu_{21} z_2 + e \mu_{31} z_3][z_1 + b\mu_{12} (2d)^{-1} z_2 
         + c\mu_{13} (2e)^{-1} z_3]\\[3mm]
    &=& b z_2^2 + c z_3^2 + 2 d z_1 z_2 + 2 e z_1 z_3 + 
         [ be\mu_{12}\mu_{23} (d\mu_{13})^{-1} + 
	         cd \mu_{13} (e\mu_{12})^{-1} ] z_2 z_3.
\end{array}
\]
\quad\\[-1mm]
If $Q$ factors, then the coefficients of $z_2^2$, $z_3^2$, $z_1 z_2$ and 
$z_1 z_3$ of $Q$ imply that $Q = Q_1$ or $Q = Q_2$. By comparing the
coefficients of $z_2 z_3$ in each case, we find $D_7(M) = 0$.
Conversely, if $D_7(M) = 0$, then $Q = Q_1$ or $Q = Q_2$, so $Q$ factors.


(b)(ii) \   
Suppose $a \ne 0$ and that $Q$ factors.  We may write 
\[
Q = a^{-1}(az_1 + \alpha_2 z_2 + \alpha_3 z_3)(az_1 + \beta_2 z_2 +
\beta_3 z_3),
\]
for some $\alpha_2$, $\alpha_3$, $\beta_2$, $\beta_3 \in \kk$.
Comparing coefficients, we have
\begin{gather}
ab = \alpha_2\beta_2,\qquad 2d = \beta_2 + \mu_{12}\alpha_2,\qquad
2e = \beta_3 + \mu_{13}\alpha_3, \tag{ix}\\
ac = \alpha_3\beta_3, \qquad
2af = \alpha_2\beta_3 + \mu_{23}\alpha_3\beta_2.\tag{x}
\end{gather}
Equations \thetag{ix} imply that $ab = \alpha_2(2d - \mu_{12}\alpha_2)$,
and so $\alpha_2 = \mu_{21}(d + X)$, where $X^2 = d^2 - \mu_{12} ab$.
Similarly, $\alpha_3 = \mu_{31}(e + Y)$, where $Y^2 = e^2 - \mu_{13} ac$.

 From the second equation in \thetag{x}, it follows that\\[-3mm]
\[
\begin{array}{rcl}
2af &=& \alpha_2(2e - \mu_{13}\alpha_3) + 
             \mu_{23}\alpha_3(2d - \mu_{12}\alpha_2)\\[3mm]
    &=& \mu_{21}(d + X)(e - Y) + \mu_{23}\mu_{31}(d - X)(e + Y),
\end{array}
\]
\quad\\[-3mm]
where $X$ and $Y$ are as above.
Hence, $D_8(M) = 0$ for some $X$ and $Y$ such that 
$X^2 = d^2 - \mu_{12} ab$ and $Y^2 = e^2 - \mu_{13} ac$.

Conversely, suppose $a \neq 0$ and that $D_8(M) = 0$ for some $X$ and $Y$ 
satisfying $X^2 = d^2 - \mu_{12} ab$ and $Y^2 = e^2 - \mu_{13} ac$. Let 
$Q' \in S_2$, where\\[-3mm] 
\[
\begin{array}{rcl}
Q' &=& a^{-1}[az_1 + \mu_{21}(d + X)z_2 + \mu_{31}(e + Y)z_3] 
       [az_1 + (d - X)z_2 + (e - Y)z_3]\\[3mm]
   &=& az_1^2 + bz_2^2 + cz_3^2 + 2dz_1z_2 + 2ez_1z_3 +\\[2mm] 
   & &\qquad\quad 
        + a^{-1}[\mu_{21}(d+X)(e-Y) + \mu_{23} \mu_{31}(e+Y)(d-X)] z_2z_3.
\end{array}
\]
\quad\\[-2mm]
The last coefficient equals $2f$, since $D_8(M) = 0$, and so $Q' = Q$,
which completes the proof of (b)(ii).
\end{pf}

We remark that, in Theorem~\ref{muRk3VarsThm}(b)(ii), converting the equation 
$D_8(M) = 0$ to a polynomial equation yields, at best, a user-unfriendly
polynomial equation of degree six:\\[-3mm]
\[
\begin{array}{l}
0 = (\mu_{13} + \mu_{12} \mu_{23})^4 a^2 b^2 c^2 
+64 \mu_{12} \mu_{13} \mu_{23} d^2 e^2 f^2 +\\[2mm]
\qquad
+ 16 (\mu_{12}^2 \mu_{13}^2 a^2 f^4 + \mu_{12}^2 \mu_{23}^2 b^2 e^4 
+ \mu_{13}^2 \mu_{23}^2 c^2 d^4)+\\[2mm]
\qquad
+ 16 (\mu_{13}^2 + \mu_{12}^2 \mu_{23}^2) (\mu_{12} a b e^2 f^2 
+ \mu_{13} a c d^2 f^2 + \mu_{23} b c d^2 e^2 )+\\[2mm]
\qquad -32 (\mu_{13} + \mu_{12} \mu_{23}) (\mu_{12} \mu_{13} a d e f^3 
+ \mu_{12} \mu_{23} b d e^3 f + \mu_{13} \mu_{23}  c d^3 e f )+\\[2mm]
\qquad
-8 (\mu_{13} + \mu_{12} \mu_{23})^2 ( \mu_{12} \mu_{13} a^2 b c f^2
+ \mu_{12} \mu_{23} a b^2 c e^2 + \mu_{13} \mu_{23} a b c^2 d^2)+\\[2mm]
\qquad
-8(\mu_{13}^3 - 5 \mu_{12} \mu_{13}^2 \mu_{23} - 5 \mu_{12}^2 \mu_{13}
            \mu_{23}^2 + \mu_{12}^3 \mu_{23}^3)a b c d e f.
\end{array}
\]
\quad\\

Theorem~\ref{muRk3VarsThm} suggests the following generalization of 
$\mu$-rank in Definition~\ref{muRk2VarsDefn} to the three-generator case.

\begin{defn}\label{muRank3Vars}
Let $Q = az_1^2 + bz_2^2 + cz_3^2 + 2dz_1z_2 + 2ez_1z_3 + 2fz_2z_3 \in
S_2$, where $a, \ldots , f \in \kk$, with $a = 0$ or 1, 
let $M \in M^{\mu}(3,\ \kk)$ be the $\mu$-symmetric matrix associated to $Q$ 
and let $D_i : M^{\mu}(3,\ \kk) \to \kk$, for $i = 1, \ldots , 8$, be 
defined as in Definition~\ref{defnDi}.
If $n=3$, we define the function $\mu$-rank $: S_2 \to \N$ as follows:
\begin{enumerate}[(a)]
  \item if $Q = 0$, we define $\mu$-rank$(Q) = 0$;
  \item if $Q \ne 0$ and if $D_i(M) = 0$ for all $i = 1, \ldots , 6$, 
  we define $\mu$-rank$(Q) = 1$;
  \item if $D_i(M) \ne 0$ for some $i = 1, \ldots , 6$ and if 
          \[(1 - a)D_7(M) + aD_8(M) = 0,\] 
	  we define $\mu$-rank$(Q) = 2$;
  \item if $(1 - a)D_7(M) + aD_8(M) \neq 0$, we define $\mu$-rank$(Q) = 3$.
\end{enumerate}
\end{defn}

\begin{example} 
If $ Q =  (2z_1 + z_2 + 8z_3)^2 = 
    (2\mu_{12}z_1 + z_2 + 8z_3)(2\mu_{21}z_1 + z_2 + 8 z_3)$,
    where $\mu_{12} = \mu_{13}$, 
    then $\mu$-rank$(Q) = 1$, 
by Definition~\ref{muRank3Vars} and Theorem~\ref{muRk3VarsThm}(a).
\end{example}

\begin{cor}\label{factoring} Let $n=3$.
\begin{enumerate}
\item[{\rm (a)}]
If $Q \in S_2^\times$, then $\mu$-rank$(Q) \leq 2$ if and only if 
$Q = L_1 L_2$ for some $L_1$, $L_2 \in S_1^\times$. 
\item[{\rm (b)}]
If $Q \in S_2^\times$, then $\mu$-rank$(Q) = 1$ if and only if $Q = L^2$ 
for some $L \in S_1^\times$. 
\end{enumerate}
\end{cor}
\begin{pf}
The result follows from Theorem~\ref{muRk3VarsThm}.
\end{pf}



\bigskip

The following result gives simplified versions of $D_7$  and $D_8$
in the special case where $S$ is a twist of the polynomial ring (see
Remark~\ref{twistrmk}).

\begin{cor}
Let $n=3$.
If $S$ is a twist of the polynomial ring by an automorphism $($see
Remark~\ref{twistrmk}$)$, then 
\[ D_7(M) = (\mu_{23}cd^2 - 2def + be^2)^2
\qquad \text{and} \qquad
D_8(M) = 2[\, \mu_{21}(de - XY )- af\, ],  \]
where $X^2 = d^2 - \mu_{12} ab$ and $Y^2 = e^2 - \mu_{13} ac$. 
\end{cor}
\begin{pf}
By Remark~\ref{twistrmk}, $\mu_{13} = \mu_{12}\mu_{23}$, so the result
follows.
\end{pf}

\bigskip

The results in this article suggest that generalizing the notion of rank to
quadratic forms on four or more generators is likely to be very computation 
heavy.  However, in the spirit of Corollary~\ref{factoring}, one could define 
$\mu$-rank one, respectively $\mu$-rank two, by simply using factoring
as follows.
\begin{defn}\label{rank1&2}
Let $n \in \N$, $n > 0$, and let $Q \in S_2^\times$.
\begin{enumerate}
\item[{\rm (a)}]
If $Q = L^2$ for some $L \in S_1^\times$, we define $\mu$-rank$(Q) = 1$.
\item[{\rm (b)}]
If $Q \neq L^2$ for any $L \in S_1^\times$, but $Q = L_1 L_2$ where 
$L_1$, $L_2 \in S_1^\times$, we define 
$\mu$-rank$(Q) = 2$.
\end{enumerate}
\end{defn}


\bigskip

\begin{thebib}{99}
\raggedbottom
\itemsep7pt
\baselineskip18pt
\bibitem[ATV1]{ATV1} 
{\sc M.~Artin,  J.~Tate and M.~Van den Bergh},  Some Algebras
Associated to Automorphisms of Elliptic Curves, {\it The Grothendieck
Festschrift} {\bf 1}, 33-85,
Eds.\ P.\ Cartier et al., Birkh\"auser (Boston, 1990).

\bibitem[ATV2]{ATV2} 
{\sc M.~Artin, J.~Tate and M.~Van den Bergh}, Modules over Regular
Algebras of Dimension~3, {\it Invent.\ Math.} {\bf 106} (1991), 335-388.

\bibitem[AL]{Aubry.Lemaire}
{\sc M.\ Aubry and J.-M.\ Lemaire}, Zero Divisors in Enveloping Algebras of
Graded Lie Algebras, {\it J.~Pure and App.\ Algebra} {\bf 38} (1985),
159-166.

\bibitem[CV]{CV} 
{\sc T.\ Cassidy and M.\ Vancliff}, Generalizations of Graded Clifford
Algebras and of Complete Intersections, {\em J.\ Lond.\ Math.\ Soc.}\ 
{\bf 81} (2010), 91-112.

\bibitem[L]{L}
{\sc L.~Le Bruyn}, Central Singularities of Quantum Spaces, {\it J.
Algebra} {\bf 177} No.~1 (1995), 142-153.

\bibitem[N]{N} 
{\sc M.\ Nafari}, {\em Regular Algebras Related to Regular Graded Skew
Clifford Algebras of Low Global Dimension}, Ph.D.\ Thesis, University of
Texas at Arlington, August 2011.

\bibitem[NVZ]{NVZ}
{\sc M.~Nafari, M.~Vancliff and Jun Zhang}, Classifying Quadratic Quantum
$\PP^2$s by using Graded Skew Clifford Algebras, {\it J.~Algebra} {\bf 346}
No.~1 (2011), 152-164.

\bibitem[VVW]{VVW} 
{\sc M.~Vancliff, K.~Van Rompay and L.~Willaert}, Some Quantum $\PP^3$s
with Finitely Many Points, {\it Comm.\ Alg.} {\bf 26} No.~4 (1998),
1193-1208.

\bibitem[VV]{VcV2}
{\sc M.~Vancliff and P.~P.~Veerapen}, Point Modules over Graded Skew
Clifford Algebras, work in progress, 2012.

\bibitem[VdB]{VdB.20pts}
{\sc M.~Van den Bergh}, An Example with 20 Points, Notes (1988).
\end{thebib}

\end{document}